\documentclass[a4paper]{amsart}

\usepackage{amsmath}
\usepackage{amsthm}
\usepackage{amssymb}
\usepackage{mathrsfs}
\usepackage[latin1]{inputenc}
\usepackage{diagrams}
\usepackage{a4wide}

\title{Deformations of glued $G_{2}$-manifolds}
\author{Johannes Nordström}

\address{Department of Mathematics, Imperial College London,
	London SW7 2AZ, United Kingdom}
\email{j.nordstrom@imperial.ac.uk}

\newcommand{\ignore}[1]{}

\DeclareMathOperator{\re}{re}

\DeclareMathOperator{\im}{im}
\DeclareMathOperator{\diam}{diam}

\newcommand{\half}{{\textstyle\frac{1}{2}}}

\newcommand{\bbr}{\mathbb{R}}

\newcommand{\bbrp}{\mathbb{R}^{+}}

\newcommand{\into}{\hookrightarrow}

\newcommand{\gtstr}{$G_{2}$\nobreakdash-\hspace{0pt}structure}

\newcommand{\gtmfd}{$G_{2}$\nobreakdash-\hspace{0pt}manifold}

\newcommand{\gtmetric}{$G_{2}$\nobreakdash-\hspace{0pt}metric}

\newcommand{\ltwoorth}{$L^{2}$\nobreakdash-\hspace{0pt}orthogonal}

\newcommand{\cystr}{Calabi-Yau structure}

\newcommand{\harm}{\mathcal{H}}
\newcommand{\norm}[1]{\Vert #1 \Vert}

\newcommand{\xcal}[1]{\mathcal{#1}_{X}}

\newcommand{\cald}{\mathcal{D}}

\newcommand{\calr}{\mathcal{R}}
\newcommand{\calq}{\mathcal{Q}}

\newcommand{\calx}{\mathcal{X}}

\newcommand{\calm}{\mathcal{M}}
\newcommand{\caln}{\mathcal{N}}

\newcommand{\calg}{\mathcal{G}}
\newcommand{\cala}{\mathcal{A}}
\newcommand{\cale}{\mathcal{E}}
\newcommand{\calb}{\mathcal{B}}

\newcommand{\defstr}{\calm}

\newcommand{\defstrbar}{\overline\defstr}
\newcommand{\calgbar}{\overline\calg}
\newcommand{\defstrgi}{\defstr_{\gi}}

\newcommand{\tv}{\tilde \varphi}
\newcommand{\tp}{\tilde \psi}

\newcommand{\diag}{d}
\newcommand{\gi}{\pm}

\newcommand{\diagcal}[1]{\mathcal{#1}_{\diag}}
\newcommand{\ycal}[1]{\mathcal{#1}_{y}}
\newcommand{\gical}[1]{\mathcal{#1}_{\gi}}
\newcommand{\onecal}[1]{\mathcal{#1}_{+}}
\newcommand{\twocal}[1]{\mathcal{#1}_{-}}
 
\newtheorem{thm}{Theorem}[section]
\newtheorem{prop}{Proposition}[section]
\newtheorem{lem}{Lemma}[section]

\theoremstyle{definition}
\newtheorem{defn}{Definition}[section]
\theoremstyle{remark}
\newtheorem{rmk}{Remark}[section]

\numberwithin{equation}{section}

\begin{document}

\begin{abstract}
We study how a gluing construction, which produces compact manifolds with
holonomy $G_{2}$ from matching pairs of asymptotically cylindrical \gtmfd s,
behaves under deformations.
We show that the gluing construction defines a smooth map
from a moduli space of gluing data to the moduli space $\defstr$
of torsion-free \gtstr s on the glued manifold,
and that this is a local diffeomorphism.
We use this to partially compactify $\defstr$, including it as the interior of
a topological manifold with boundary.
The boundary points are equivalence classes of matching pairs of torsion-free
asymptotically cylindrical \gtstr s.
\end{abstract}

\maketitle

\section{Introduction}

The exceptional Lie group $G_{2} \subset SO(7)$ also occurs as an exceptional
case in the classification of Riemannian holonomy groups due to Berger
\cite{berger55}.  A \gtmfd{} is a $7$-dimensional Riemannian manifold with
holonomy group contained in~$G_{2}$. Its metric can be defined in terms of a
closed differential \mbox{$3$-form} equivalent to a torsion-free \gtstr. Joyce
\cite{joyce96-I} constructed the first compact examples of manifolds with
holonomy $G_{2}$. He also proved that the moduli space $\defstr$ of
torsion-free \gtstr s on a compact \gtmfd, i.e. the quotient of the space of
torsion-free \gtstr s by the identity component of the diffeomorphism group, is
a smooth manifold.

A \gtmfd{} is asymptotically cylindrical if it is
asymptotically isometric to a product cylinder outside a compact subset.
Kovalev \cite{kovalev03} explains a gluing construction which produces a
compact \gtmfd{} $M$ from a pair $M_{\gi}$ of asymptotically cylindrical
\gtmfd s with matching cylindrical parts, and constructs
new examples of compact manifolds with holonomy $G_{2}$ by this method.
Topologically, $M$ can be considered
as a generalised connected sum of $M_{+}$ and~$M_{-}$.
The author \cite{jn1} shows that there is a smooth moduli space of torsion-free
\gtstr s on asymptotically cylindrical \gtmfd s, extending the result of Joyce
from the compact case. This leads to the question of how deformations of a
compact \gtmfd{} constructed by gluing are related to deformations of the
asymptotically cylindrical halves. We find that the torsion-free \gtstr s
obtainable from the gluing construction form an open subset of the moduli space
$\defstr$ on the compact manifold. This subset can be regarded as a
neighbourhood of a boundary component for~$\defstr$. 

The data required for the gluing construction is a pair
$(\varphi_{+}, \varphi_{-})$ of asymptotically cylindrical \gtstr s on $M_{+}$
and $M_{-}$ which satisfies a matching condition
(cf. definition \ref{matchdef}), together with a large
parameter $L \in \bbrp$. $L$~controls the length of an approximately
cylindrical neck in the result of the gluing. Given such a set of gluing data 
$(\varphi_{+}, \varphi_{-}, L)$, the gluing construction yields a torsion-free
\gtstr{} $Y(\varphi_{+}, \varphi_{-}, L)$ on the connected sum $M$.
This is uniquely defined up to diffeomorphisms, and so represents a well-defined
point in the moduli space $\defstr$ of torsion-free \gtstr s on $M$.
We will study the local properties of a gluing map defined on a quotient
$\calg$ of a space of gluing data by a natural symmetry group.
By relating $\calg$ to the moduli spaces of torsion-free \gtstr s on the halves
$M_{\gi}$, which are smooth manifolds, we show that $\calg$ is smooth too.
The main theorem \ref{maingluethm} states that
\[ Y : \calg \to \calm \]
is a local diffeomorphism. This result can be compared with 
Donaldson and Kronheimer's description \cite[\S 7.2]{donaldson97} of
deformations of anti-self-dual connections on a connected sum of a pair of
smooth $4$-manifolds.

We also explain how to apply these results
to attach a boundary to~$\defstr$, forming a topological manifold $\defstrbar$
with boundary, so that the boundary points correspond to ways
of `pulling apart' $M$ into a pair of asymptotically cylindrical \gtmfd s.
The results about the gluing map can therefore be interpreted as a description
of a neighbourhood of a boundary component of~$\defstr$.
Like the statement that $\defstr$ is a
manifold, this is essentially a local result. Little is known about
the global properties of~$\defstr$. Its local properties are
also studied for instance by Karigiannis and Leung \cite{karigiannis07} and
Grigorian and Yau \cite{grigorian08}, who consider in particular the
curvature of a natural pseudo-Riemannian metric on~$\defstr$.

The topological quantum field theory proposed by Leung \cite{leung02}
considers generalised connected sums of almost \gtmfd s, i.e. $7$-manifolds
with \gtstr{} which is not necessarily torsion-free (so the associated metric
need not have holonomy in~$G_2$) but whose defining $3$-form is still required
to be closed.
The proposed TQFT assigns invariants to compact and asymptotically cylindrical
almost \gtmfd s by counting coassociative cycles, and these invariants are
expected to behave well under connected sums. It is clear that a perturbation
of a connected sum of asymptotically cylindrical almost \gtmfd s remains such a
connected sum, but our result shows that this holds also when working in the
category of torsion-free \gtmfd s, where the metric has holonomy in~$G_2$.

The paper is organised as follows.
Section \ref{gluesetupsec} contains background for the gluing construction
of compact \gtmfd s and precise statements of our main results.
In section \ref{gluetoposec}, we discuss
the topology of the glued manifold $M$ and prove theorem \ref{gluehodgethm},
a Hodge theory gluing result of some potential independent interest.
This is used in section \ref{glueproofsec} to compute the
derivative of the gluing map, proving the main theorem \ref{maingluethm}.
In section \ref{glueboundarysec}, we outline how to attach boundary points
to $\defstr$.

\section{Setup}
\label{gluesetupsec}

\subsection{Preliminaries}

We review the preliminary definitions that are required to set up the gluing
construction and state the main results.
For more detailed background on \gtmfd s see Joyce \cite[Chapter 10]{joyce00}
or Salamon \cite[Chapter 8]{salamon89}.

Recall that $G_{2}$ can be defined as the automorphism group of the normed
algebra of octonions. Equivalently, $G_{2}$ is the stabiliser in
$GL(\bbr^{7})$ of the $3$-form
\begin{equation}
\label{g2formeq}
\varphi_{0} = dx^{123} + dx^{145} + dx^{167} + dx^{246}
 - dx^{257} - dx^{347} - dx^{356} \in \Lambda^{3}(\bbr^{7})^{*} .
\end{equation}
A \emph{\gtstr{}} on a manifold $M^{7}$ can therefore be defined in terms of a
differential $3$-form $\varphi$ which is equivalent to $\varphi_{0}$ at
each point.
$\varphi_{0}$ is a stable form in the sense of Hitchin \cite{hitchin01},
i.e. the $GL(\bbr^{7})$-orbit of $\varphi_{0}$ is open in
$\Lambda^{3}(\bbr^{7})^{*}$, so the set of \gtstr s on $M$ is open in
the space of $3$-forms on $M$ with respect to the uniform norm.
A \gtstr{} naturally defines a Riemannian metric $g_{\varphi}$
and an orientation on $M$, and thus also a Levi-Civita connection
$\nabla_{\varphi}$ and a Hodge star~$*_{\varphi}$.
$\varphi$~is called \emph{torsion-free} if $\nabla_{\varphi}\varphi = 0$.
By a result of Gray this condition is equivalent to
$d\varphi = d^{*}_{\varphi}\varphi = 0$. Note that this is a non-linear
condition, since $d^{*}_{\varphi}$ depends on $\varphi$.
We call a $7$-dimensional manifold $M$ equipped with a torsion-free \gtstr{}
and the induced Riemannian metric a \emph{\gtmfd}.

The \emph{holonomy} group of a Riemannian manifold is the group of isometries
of a tangent plane generated by parallel transport around closed curves.
Parallel tensor fields on the manifold correspond to invariants of the
holo\-nomy group, so it is clear that a $7$-dimensional Riemannian manifold $M$
has holo\-nomy $Hol(M)$ contained in $G_{2}$ if and only if the metric is
induced by a torsion-free \gtstr. For a compact \gtmfd{} the holo\-nomy is
exactly $G_{2}$ if and only if the fundamental group $\pi_{1}M$ is finite (see
\cite[Proposition 10.2.2]{joyce00}), otherwise a finite cover of $M$ is
a Riemannian product of lower-dimensional manifolds.

On a compact \gtmfd{} $M$ the group $\cald$ of diffeomorphisms isotopic to the
identity acts on the space $\calx$ of torsion-free \gtstr s by pull-backs.
The quotient $\defstr = \calx/\cald$ is the moduli space of torsion-free
\gtstr s. Since torsion-free \gtstr s are closed forms there is a natural
projection $\defstr \to H^{3}(M)$ to de Rham cohomology.

\begin{thm}[{\cite[Theorem $10.4.4$]{joyce00}}]
\label{maincptthm}
Let $M$ be a compact \gtmfd. Then the moduli space $\defstr$ of torsion-free
\gtstr s on $M$ is a smooth manifold, and the map $\defstr \to H^{3}(M)$
is a local diffeomorphism.
\end{thm}

For $X^{6}$ compact, we call a \gtstr{} on $X \times \bbr$ \emph{cylindrical}
if it is translation-invariant and defines a product metric.
The stabiliser in $G_{2}$ of a vector in $\bbr^{7}$ is $SU(3)$.
The product of a Riemannian manifold $X^{6}$ with $\bbr$ therefore
has $Hol(X \times \bbr) \subseteq G_{2}$ if and only if
$Hol(X) \subseteq SU(3)$, so the cross-section of a \emph{cylindrical \gtmfd}
is always a Calabi-Yau $3$-fold.
If we let $z^{1} = x^{2}+ix^{3}, z^{2} = x^{4}+ix^{5},z^{3} = x^{6}+ix^{7}$
then we can write $\varphi_{0}$ as
\begin{equation}
\label{pointmodeleq}
\varphi_{0} = \Omega_{0} + dx^{1} \wedge \omega_{0} ,
\end{equation}
where
\begin{gather}
\label{omega01eq}
\Omega_{0} = \re (dz^{1} \wedge dz^{2} \wedge dz^{3}) , \\
\omega_{0} = {\textstyle\frac{i}{2}}(dz^{1} \wedge d \bar z^{1} +
dz^{2} \wedge d \bar z^{2} + dz^{3} \wedge d \bar z^{3}) .
\end{gather}
A cylindrical \gtstr{} $\varphi$ on $X \times \bbr$ is therefore of the form
\[ \varphi = \Omega + dt \wedge \omega , \]
where $(\Omega, \omega)$ is a pair of forms on $X$ point-wise equivalent to 
$(\Omega_{0}, \omega_{0})$. If $\varphi$ is torsion-free then
$(\Omega, \omega)$ can be considered to define a Calabi-Yau structure
(or torsion-free $SU(3)$-structure) on $X$. This means that $X$ has an
integrable complex structure, $\omega$ is the Kähler form of a Ricci-flat
Kähler metric, and $\Omega$ is the real part of a non-vanishing holomorphic
$(3,0)$-form.

A non-compact manifold $M$ is said to have \emph{cylindrical ends}
if $M$ is written as union of two pieces $M_{0}$ and $M_{\infty}$ with
common boundary $X$, where $M_{0}$ is compact,
and $M_{\infty}$ is identified with $X \times \bbrp$
by a diffeomorphism (identifying $\partial M_{\infty}$ with $X \times \{0\}$).
$X$~is called the \emph{cross-section} of~$M$.
Let $t$ be a smooth real function on $M$ which is the $\bbrp$-coordinate
on $M_{\infty}$, and negative on the interior of $M_{0}$.
A tensor field $s$ on $M$ is said to be \emph{exponentially asymptotic}
with rate $\delta > 0$ to a translation-invariant tensor field $s_{\infty}$
on $M$ if $e^{\delta t} \norm{\nabla^{k}(s-s_{\infty})}$ is bounded
on $M_{\infty}$ for all $k \geq 0$,
with respect to a norm defined by an arbitrary Riemannian metric on $X$.

A metric on $M$ is called \emph{EAC} (exponentially asymptotically cylindrical)
if it is exponentially asymptotic to
a product metric on $X \times \bbr$, and a \gtstr{} is said to be EAC if it
is exponentially asymptotic to a cylindrical \gtstr{} on $X \times \bbr$.
The asymptotic limit of a torsion-free EAC \gtstr{} then defines a
\cystr{} on the cross-section $X$.
A diffeomorphism $\phi$ of $M$ is called EAC if it is exponentially
close to a product diffeomorphism $(x,t) \mapsto (\Xi(x), t + h)$ of
$X \times \bbr$ in a similar sense.

The moduli space of torsion-free EAC \gtstr s on an EAC \gtmfd{} $M$ is the
quotient of the space of torsion-free EAC \gtstr s (with any exponential rate)
by the group of EAC diffeomorphisms of $M$. We will review some properties
of the EAC moduli space in subsection \ref{eacdefsub}, but note for now that
theorem \ref{maincptthm} from the compact case can be generalised.

\begin{thm}[{\cite[Theorem 3.2]{jn1}}]
\label{maineacthm}
Let $M$ be an EAC \gtmfd. Then the moduli space of torsion-free EAC \gtstr s
on $M$ is a smooth manifold.
\end{thm}

\subsection{Gluing construction}
\label{glueconstrsub}

Let $M_{\pm}$ be a pair of oriented dimension $7$ manifolds, each with a single
cylindrical end, and the same cross-section $X$.
We assume that $X$ is oriented so that its orientation
agrees with that defined by $M_{+}$ on its boundary, and is the reverse of that
defined by $M_{-}$ on its boundary. This ensures that the connected sum of
$M_{+}$ and $M_{-}$ obtained by identifying their boundaries at infinity is
oriented.
Let $t_{\pm}$ be cylindrical coordinates on $M_{\pm}$ respectively.

\begin{defn}
\label{matchdef}
Let $\varphi_{\gi}$ be torsion-free EAC \gtstr s on $M_{\gi}$. The pair
$(\varphi_{+},\varphi_{-})$ is said to \emph{match} if their asymptotic models
are $\Omega \pm dt_{\pm} \wedge \omega$,
respectively, for some Calabi-Yau structure $(\Omega,\omega)$ on $X$
compatible with the chosen orientation.
Let $\ycal{X}$ be the space of such pairs.
\end{defn}

Given $L \in \bbrp$
let $M_{\gi}(L) = \{ y \in M_{\gi} : t_{\gi} \leq L \}$.
Identify the boundaries of $M_{\gi}(L)$ to form a
compact smooth manifold $M(L)$, and
let $j^{*} : X \into M(L)$ be the inclusion of the common boundary.
$M(L)$ is independent of $L$ up to diffeomorphism, so we will
often refer to it simply as~$M$.

For notational convenience we suppose that the cylindrical end of
$M_{\gi}$ is given by $t_{\gi} > -2$.
Let $\rho_{\gi}$ be a smooth cut-off function on $M_{\gi}$ which is $0$ for
$t_{\gi} < L-2$ and $1$ for $t_{\gi} > L-1$.
Let $\alpha$ be a closed exponentially asymptotically
translation-invariant $m$-form on~$M_{\gi}$. Then it can be written
as $\alpha_{\infty} + \beta_{t_\gi} + dt_\gi \wedge \gamma_{t_\gi}$ on the
cylinder, with $\alpha_{\infty}$ translation invariant, and
$\beta _{t_{\gi}} \in \Omega^{m}(X)$, $\gamma_{t_{\gi}} \in \Omega^{m-1}(X)$
both exponentially decaying in $t_{\gi}$.
Define an $(m-1)$-form on the cylinder by
\begin{equation}
\label{etaeq}
\eta_{\gi}(\alpha) = \rho_{\gi} \int_{t_{\gi}}^{\infty} \gamma_{s} ds .
\end{equation}
Then $\alpha + d\eta_{\gi}(\alpha)$ is translation-invariant on $t > L-1$.

For $(\varphi_{+},\varphi_{-}) \in \ycal{X}$
let $\tv_{\gi} = \varphi_{\gi} + d\eta_{\gi}(\varphi_{\gi})$.
Then we can define a \gtstr{} $\tv(\varphi_{+},\varphi_{-}, L)$ on $M(L)$ by
$\tv|_{M_{\gi}(L)} = \tv_{\gi}|_{M_{\gi}(L)} $.
Note that the choice of cut-off function in the definition of
$\eta_{\gi}$ does not affect the cohomology class of
$\tv(\varphi_{+},\varphi_{-}, L)$.

\begin{prop}
\label{basicglueprop}
There is an upper semi-continuous map $L_{0} : \ycal{X} \to \bbrp$ such
that for any $L > L_{0}$ there is a unique diffeomorphism class of
torsion-free \gtstr s on $M(L)$ in a small neighbourhood
of $\tv(\varphi_{+},\varphi_{-}, L)$ in its cohomology class.
\end{prop}

\begin{proof}[Sketch proof.]
The idea is that for large $L$ the torsion of $\tv(\varphi_{+},\varphi_{-}, L)$
is very small, and the structure can be perturbed to a torsion-free one using
a contraction-mapping argument. See Kovalev \cite[\S 5]{kovalev03} for details.
The argument is inspired by a construction of Floer \cite{floer91}
(see also Kovalev and Singer \cite{kovalev01}).
\end{proof}

The resulting \gtmetric{} on $M(L)$ has an almost cylindrical `neck' of length
roughly $2L$, and $\diam M(L) \sim 2L$ as $L \to \infty$.

Kovalev \cite{kovalev03} constructs examples of matching pairs of EAC \gtmfd s
to which the gluing construction can be applied. An EAC version of the
Calabi conjecture produces EAC manifolds with holonomy $SU(3)$. These can be
multiplied by circles $S^{1}$ to form (reducible) \gtmfd s, which form
compact irreducible \gtmfd s (manifolds with holonomy exactly~$G_{2}$)
when glued together. These have different topological type from the examples
constructed earlier by Joyce \cite{joyce96-I}.

A future paper \cite{jn2} will explain how some of the examples of compact
\gtmfd s constructed by Joyce can also be produced by gluing a pair of
EAC \gtmfd s. In some of these examples the EAC components are
irreducible EAC \gtmfd s.

\subsection{Statement of results}

Let $M$ be the gluing of two EAC \gtmfd s $M_{\gi}$ as above.
Let $\defstr$ be the moduli space of torsion-free \gtstr s on $M$, and
$\defstr_{\gi}$ the moduli spaces of torsion-free EAC \gtstr s on $M_{\gi}$.
These are all smooth manifolds by theorems \ref{maincptthm} and
\ref{maineacthm}.

When considering how the gluing construction behaves under deformations
it is natural to look at the space of matching pairs of diffeomorphism
classes of torsion-free \gtstr s on $M_{+}$ and~$M_{-}$, i.e.
the subset $\ycal{M} \subseteq \onecal{M} \times \twocal{M}$
consisting of pairs which have matching images in the moduli space
of Calabi-Yau structures on $X$. We will
use the deformation theory for EAC \gtmfd s from \cite{jn1} to show that
$\ycal{M}$ is a manifold. However, given a matching pair of diffeomorphism
classes of EAC \gtstr s there is some ambiguity in how to glue them, since we
need to choose how to identify the cylindrical ends.
This means both choosing how to identify the
cross-sections (this ambiguity roughly corresponds to the quotient of
the automorphism group of the cross-section by a subgroup generated by
elements which extend to automorphisms of $M_{+}$ or~$M_{-}$), and
choosing the neck length for the glued manifold.
It is therefore not possible to use $\ycal{M}$ itself as the domain for
any sensible, single-valued map to $\defstr$.
Instead we define a gluing map on a moduli space of data for the gluing
construction.

\begin{defn}
A set of \emph{gluing data} is a triple
$(\varphi_{+},\varphi_{-},L) \in \ycal{X} \times \bbr$ such that 
$(\varphi_{+},\varphi_{-}) \in \ycal{X}$ and
$L > L_{0}(\varphi_{+},\varphi_{-})$.
Let $G_{0}$ be the space of gluing data.
\end{defn}

$G_{0}$ is an open subset of $\ycal{X} \times \bbr$.
Proposition \ref{basicglueprop} provides a well-defined smooth
map to the moduli space of torsion-free \gtstr s on $M$,
\begin{equation}
\label{upperyeq}
Y : G_{0} \to \mathcal{M} .
\end{equation}

Two sets of gluing data define essentially the same gluing operation if they
are equivalent under the following action. Let $\gical{D}$ be the group of EAC
diffeomorphisms of $M_{\gi}$ isotopic to the identity.

\begin{defn}
\label{calyddef}
$(\phi_{+}, \phi_{-}) \in \onecal{D} \times \twocal{D}$ such that
$\phi_{\gi}$ is asymptotic to
$(x,t_{\gi}) \mapsto (\Xi_{\gi}(x), t_{\gi} + h_{\gi})$
is a \emph{matching pair of EAC diffeomorphisms} if
$\Xi_{+} = \Xi_{-}$. Let $\ycal{D}$ be the identity component
of the group of such pairs.

For  $(\phi_{+}, \phi_{-}) \in \ycal{D}$ let $h = \half(h_{+}+h_{-})$, and
define an action on $\ycal{X} \times \bbr$ by
\begin{equation}
\label{calydacteq}
\phi^{*} : (\varphi_{+},\varphi_{-},L) \mapsto
(\phi_{+}^{*}\varphi_{+}, \phi_{-}^{*}\varphi_{-}, L-h) .
\end{equation}
\end{defn}
There is no reason why the open set
$G_{0} \subseteq \ycal{X} \times \bbr$
should be invariant under the action of~$\ycal{D}$. Nevertheless we can define

\begin{defn}
The \emph{moduli space of gluing data} is
$\mathcal{G}_{0} = G_{0}\ycal{D}/\ycal{D}$.
\end{defn}

We can project (\ref{calydacteq}) to an action of $\ycal{D}$ on
the space $\ycal{X}$ of matching pairs. The quotient $\calb$ has a natural map
to $\ycal{M}$. By studying this map we will deduce smoothness of $\calb$ from
the fact that $\ycal{M}$ is a smooth manifold.
$\calg_{0}$~is obviously a fibre bundle over $\calb$ with typical
fibre $\bbrp$, so it is a smooth manifold too.
Moreover, the gluing map
(\ref{upperyeq}) really is invariant under the action of~$\ycal{D}$, and
therefore descends to a smooth map
\begin{equation}
\label{mainyeq}
Y : \mathcal{G}_{0} \to \mathcal{M} .
\end{equation}

Proposition \ref{yhderprop} computes the derivative of the gluing map
(\ref{mainyeq}). For each matching pair
$(\varphi_{+},\varphi_{-})$ the derivative is invertible at
$(\varphi_{+},\varphi_{-}, L)\ycal{D} \in \calg_{0}$ for all
large $L$. Therefore $Y$ is a local diffeomorphism on
some open subset $\mathcal{G} \subseteq \mathcal{G}_{0}$ whose gluing
parameters are sufficiently large.
This gives our main result.

\begin{thm}
\label{maingluethm}
Let $M$ be a compact \gtmfd{} constructed by gluing a matching pair
$M_{\gi}$ of EAC \gtmfd s. Then the gluing space $\calg$ is a smooth manifold,
and the gluing map $Y : \mathcal{G} \to \mathcal{M}$ is a local diffeomorphism.
\end{thm}

In the proof we will assume that $b^{1}(M) = 0$ in order to simplify some
technical statements; for example, the map $\calb \to \ycal{M}$ is then a
covering map.  This is not a very restrictive assumption, since $b^{1}(M) = 0$
when $M$ has holonomy exactly $G_{2}$, which is the most interesting case.  In
general $\calb \to \ycal{M}$ is a submersion, and the fibres have dimension
$b^{1}(M)$.

The most important tool in the proof is to use the local diffeomorphism
$\pi_{H} : \mathcal{M} \to H^{3}(M)$. This means that we can study the local
properties of the gluing map in terms of what the gluing does to the cohomology
classes. This is discussed in section \ref{gluetoposec};
in particular we prove a Hodge theory gluing result.

Theorem \ref{maingluethm} is proved in section \ref{glueproofsec}.
In section \ref{glueboundarysec} we outline how these arguments
can also be used to show that $\defstr$ can be partially
compactified by inclusion in a topological manifold $\defstrbar$ with boundary,
so that the paths defined by gluing a matching pair of EAC \gtstr s
with increasing gluing parameter converge to a boundary point.
The boundary points can therefore be considered as ways of `pulling apart' $M$
into a pair of EAC connected-summands.
Since the gluing space $\calg$ is a fibre bundle over $\calb$ with typical
fibre $\bbrp$, there is a natural way to
form a fibre bundle $\calgbar$ over $\calb$ with typical fibre $(0,\infty]$,
and $\partial\calgbar = \calb$.
The partial compactification of $\defstr$ can then be described in the
following way.

\begin{thm}
\label{gluebdrythm}
Let $M$ be a compact \gtmfd{} constructed by gluing a matching pair
$M_{\gi}$ of EAC \gtmfd s. Then the moduli space $\defstr$ of torsion-free
\gtstr s on $M$ can be included as the interior of a topological manifold
$\defstrbar$ with a boundary $\partial\defstrbar$, so that the gluing map
$Y$ extends to a local homeomorphism
\[ Y : \calgbar \to \defstrbar . \]
The restriction of $Y$ to the boundary is a covering map
$\partial\calgbar \to \partial\defstrbar$.
\end{thm}

\section{Gluing and topology}
\label{gluetoposec}

\subsection{Topology of the connected sum}

Let $M^{n}_{+}, M^{n}_{-}$ be oriented manifolds, each with a single
cylindrical end, which have common cross-section $X^{n-1}$. As in
subsection \ref{glueconstrsub} we assume that $X$ is oriented compatible with
$M_{+}$ and reverse to $M_{-}$, and we form a generalised connected sum $M$.
We collect here some results about the topology of $M$ that we will use.

As we remarked before, as a smooth manifold $M$ is independent of the choice of
gluing para\-meter~$L$.
Up to isotopy there are natural inclusion maps

%\begin{diagram}[nohug,PostScript=Rikicki]
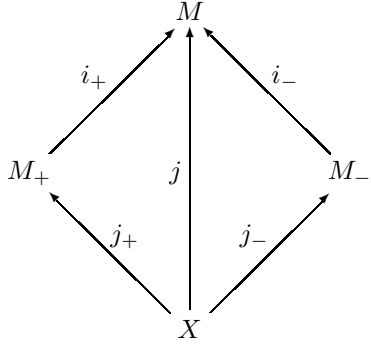
\begin{diagram}[nohug,notextflow]
& & M & & \\
& \ruTo^{i_{+}} & & \luTo^{i_{-}} & \\
M_{+} & & \uTo^{j} & & M_{-} \\
& \luTo^{j_{+}} & & \ruTo^{j_{-}} & \\
& & X & &
\end{diagram}

A large part of what we need to know about the topology is contained in the
exactness of the Mayer-Vietoris sequence for $M = M_{+} \cup M_{-}$ and the
sequence for the cohomology of $M_{\gi}$ relative to its boundary $X$.
Throughout, $H^{*}$ refers to de Rham cohomology.
\begin{gather}
\label{mayereq}
\cdots \longrightarrow H^{m-1}(X) \stackrel{\delta}{\longrightarrow}
H^{m}(M) \stackrel{i_{+}^{*} \oplus i_{-}^{*}}{\longrightarrow}
H^{m}(M_{+}) \oplus H^{m}(M_{-})
\stackrel{j_{+}^{*} - j_{-}^{*}}{\longrightarrow} H^{m}(X) \longrightarrow
\cdots \\
\label{relexacteq2}
\cdots \longrightarrow H^{m-1}(X) \stackrel{\partial_{\gi}}{\longrightarrow}
H^{m}_{cpt}(M_{\gi}) \stackrel{e_{\gi}}{\longrightarrow} H^{m}(M_{\gi})
\stackrel{j_{\gi}^{*}}{\longrightarrow} H^{m}(X) \longrightarrow \cdots
\end{gather}
Here $e_{\gi}$ is induced by the natural chain map
$\Omega^{*}_{cpt}(M_\gi) \to \Omega^{*}(M_\gi)$, and $\delta$ and
$\partial_{\gi}$ denote the boundary homomorphisms.
The inclusions $i_{\gi} : M_{\gi} \hookrightarrow M$ induce maps
$i_{\gi*} : H^{m}_{cpt}(M_{\gi}) \to H^{m}(M)$. Note that
\begin{equation}
\label{dieq}
\delta = i_{+*} \circ \partial_{+} = - i_{-*} \circ \partial_{-} .
\end{equation}
$j_{\gi}^{*} : H^{m}(M_{\gi}) \to H^{m}(X)$ is the Poincar\'e dual of
$\pm \partial_{\gi} : H^{n-m-1}(X) \to H^{n-m}(M)$
(the sign difference comes from our assumption on the orientations of
$M_{\gi}$ and $X$).
The Poincar\'e dual of the Mayer-Vietoris sequence is the sequence for
the relative cohomology of $(M,X)$,
\begin{equation}
\cdots \; H^{m-1}(X)
\stackrel{\partial_{+} \oplus \partial_{-}}{\longrightarrow}
H^{m}_{cpt}(M_{+}) \oplus H^{m}_{cpt}(M_{-})
\stackrel{i_{+*} + i_{-*}}{\longrightarrow} H^{m}(M)
\stackrel{j^{*}}{\longrightarrow} H^{m}(X) \; \cdots
\end{equation}

Denote the image of $j_{\gi}^{*} : H^{m}(M_{\gi}) \to H^{m}(X)$ by
$A^{m}_{\gi}$, and let $A^{m}_{\diag}$ be the image of
$j^{*} : H^{m}(M) \to H^{m}(X)$. By the exactness of the Mayer-Vietoris
sequence, $A^{m}_{\diag} =  A^{m}_{+} \cap A^{m}_{-}$.

\subsection{Gluing and cohomology}
\label{gluecohsub}

We explain how to glue a matching pair of closed forms on $M_{+}, M_{-}$
to a well-defined cohomology class on $M$.

Let $\ycal{Z}^{m}$ be the space of matching pairs of closed
exponentially asymptotically translation-invariant $m$-forms
on $M_{+}, M_{-}$, i.e. $(\psi_{+}, \psi_{-})$ such that $\psi_{\gi}$
is a closed exponentially asymptotically translation-invariant
$m$-form on~$M_{\gi}$,
with asymptotic limits $B_a(\psi) \pm dt_\gi \wedge B_e(\psi)$ respectively.

If $(\psi_{+}, \psi_{-}) \in \ycal{Z}^{m}$ and $L > 0$ let
$\tp_{\gi} = \psi_{\gi} + d\eta_{\gi}(\psi_{\gi})$.
Choose the cut-off function for the cylinders in the definition (\ref{etaeq})
of $\eta_{\gi}$ to ensure that $\tp_{\gi}$ is translation-invariant
on $t_{\gi} > 0$.
Then we can define $\tp(\psi_{+}, \psi_{-}, L)$ on $M(L)$ by $i^{*}_{\gi}\tp = \tp_{\gi}$.
We define a gluing map
\begin{equation}
\label{yhdefeq}
Y_{H} : \ycal{Z}^{m} \times \bbrp \to H^{m}(M), \;\:
(\psi_{+}, \psi_{-}, L) \mapsto [\tp] .
\end{equation}
$Y_{H}$ is independent of the choice of $\eta_{\gi}$ and hence well-defined.
Furthermore, we find that
$Y_{H}$ is invariant under the action of the group $\ycal{D}$
of matching diffeomorphisms from definition \ref{calyddef}.

\begin{defn}
For $(\phi_{+}, \phi_{-}) \in \ycal{D}$ with asymptotic models
$(x, t_{\gi}) \mapsto (\Xi(x), t_{\gi} + h_{\gi})$
let $h = \half(h_{+}+h_{-})$, and
define an action on $\ycal{Z}^{m} \times \bbr$ by
\begin{equation}
\label{calydacteq2}
\phi^{*} : (\psi_{+},\psi_{-},L) \mapsto
(\phi_{+}^{*}\psi_{+}, \phi_{-}^{*}\psi_{-}, L-h) .
\end{equation}
\end{defn}

\begin{prop}
\label{yhinvprop}
If $(\psi_{+}, \psi_{-}, L) \in \ycal{Z}^{m} \times \bbrp$, and
$(\phi_{+}, \phi_{-}) \in \ycal{D}$ with $h_{\gi} < L$ then
\[ Y_{H}(\psi_{+}, \psi_{-}, L) = 
Y_{H}(\phi_{+}^{*}\psi_{+}, \phi_{-}^{*}\psi_{-}, L-h)
\in H^{m}(M) . \]
\end{prop}

\begin{proof}[Sketch proof.]
Let $\tp = \tp(\psi_{+}, \psi_{-}, L)$ and $\tp' =
\tp(\phi_{+}^{*}\psi_{+}, \phi_{-}^{*}\psi_{-}, L-h)$.
$\phi_{+}$ and $\phi_{-}$ can be approximately glued to
a diffeomorphism $\tilde\phi : M(L-h) \to M(L)$ which
pulls back $[\tp]$ to $[\tp']$.
\end{proof}

\begin{prop}
\label{lderprop}
If $(\psi_{+}, \psi_{-}) \in \ycal{Z}^{m}$ with
$B_{e}(\psi) = \tau$, $L, h \in \bbrp$ then
\begin{equation}
\label{yhnegdefeq}
Y_{H}(\psi_{+}, \psi_{-}, L + h) =
Y_{H}(\psi_{+}, \psi_{-}, L) + 2h\delta([\tau]) ,
\end{equation}
where $\delta$ is the boundary homomorphism appearing in the Mayer-Vietoris
sequence \textup{(\ref{mayereq})}.
\end{prop}

\begin{proof}
It suffices to prove the result separately for the cases when $B_{a}(\psi) = 0$
and $B_{e}(\psi) = 0$.

If $B_{e}(\psi) = 0$ pick a diffeomorphism $f : (0, L) \to (0,L+h)$ which is
$id$ on $(0,1)$ and $id + h$ on $(L-1,L)$. We can define a diffeomorphism
$M(L) \to M(L+h)$ which is the identity on the images of the compact
pieces of $M_{+}$ and $M_{-}$ in $M(L)$ and
$(x,t) \mapsto (x, f(t))$ on the cylindrical part. This pulls back
$\tp(\psi_{+}, \psi_{-}, L) \mapsto \tp(\psi_{+}, \psi_{-}, L+h)$.

If $B_{a}(\psi) = 0$ let
$c_{\gi} = \pm \tp_{\gi} - d(\rho_{\gi} t_{\gi} \tau)$,
with $\rho_{\gi}$ a cut-off function chosen so that $c_{\gi}$ has support
contained in $t_{\gi} < 1$.
By definition of the Mayer-Vietoris boundary map $\delta$, the 
form on $M(L)$ obtained by gluing $d(\rho_{+}t_{+}\tau)$ and
$-d(\rho_{-}t_{-}\tau)$ is cohomologous to
$\delta((t_{+}+t_{-})[\tau]) = 2L\delta([\tau])$ for any $L$. Hence
for any~$L$
\begin{equation}
\label{yheq}
Y_{H}(\psi_{+}, \psi_{-}, L) = i_{+*}([c_{+}]) + i_{-*}(-[c_{-}]) +
2L\delta([\tau]) .
\end{equation}
Since $i_{\gi*} : H^{m}_{cpt}(M_{\gi}) \to H^{m}(M)$ and
$\delta : H^{m-1}(X) \to H^{m}(M)$ are independent of $L$ the result follows.
\end{proof}

It is convenient to use proposition \ref{lderprop} to extend
$Y_{H}$ to negative gluing parameters in a well-defined way.

\begin{defn}
\label{yhdefn}
Define
\[ Y_{H} : \ycal{Z}^{m} \times \bbr \to H^{m}(M) \]
as (\ref{yhdefeq}) on $\ycal{Z} \times \bbrp$, and extend
for any $L > 0$ and $h \in \bbr$ by (\ref{yhnegdefeq}).
\end{defn}

\subsection{EAC Hodge theory}

Let $M_{\gi}$ be an EAC manifold with cross-section $X$.
We summarise the Hodge theory from \cite[\S $5$]{jn1}
(see also Melrose \cite[\S 6.4]{melrose94}).

Let $\harm^{m}_{\gi,0}$ be the bounded harmonic forms on $M_{\gi}$.
This is a finite-dimensional space, and its elements are smooth, closed,
coclosed, and exponentially asymptotically translation-invariant.
The asymptotic limit of $\beta \in \harm^{m}_{\gi,0}$ is
a translation-invariant harmonic form
on the cylinder $X \times \bbr$, so if $\harm^{m}_{X}$ denotes the space of
harmonic $m$-forms on $X$ then the limit can be written as
\[ B_{\gi}(\beta) = B_{\gi, a}(\beta) + dt_{\gi} \wedge B_{\gi, e}(\beta)
\in \harm^{m}_{X} + dt_{\gi} \wedge \harm^{m-1}_{X} . \]
Note that
\[ j_{\gi}^{*}[\beta] = [B_{\gi,a}(\beta)] \in H^{m}(X) . \]
The image of $B_{\gi, a} : \harm^{m}_{\gi,0} \to \harm^{m}_{X}$ is therefore
precisely the space of harmonic representatives $\cala^{m}_{\gi}$ of
the cohomology classes $A^{m}_{\gi} \subseteq H^{m}(X)$.
Let $\harm^{m}_{\gi,abs} = \ker B_{\gi, e}$, and $\harm^{m}_{E}$ the subset of
exact forms in $\harm^{m}_{\gi,0}$. Then the natural map
\[ \harm^{m}_{\gi,abs} \to H^{m}(M) \]
is an isomorphism, and
\[ \harm^{m}_{\gi,0} = \harm^{m}_{\gi,abs} \oplus \harm^{m}_{\gi,E} . \]
$B_{\gi,e}$ maps $\harm^{m}_{\gi,E}$ isomorphically to its image
$\cale^{m}_{\gi}$. Further
\[ \harm^{m}_{X} = \cala^{m}_{\gi} \oplus \cale^{m}_{\gi} , \]
and this direct sum is orthogonal.

\subsection{Hodge theory and gluing}
\label{gluehodgesub}

Now suppose that $M_{\gi}$ are EAC Riemannian manifolds whose
cylindrical models match. We wish to consider what the
gluing of closed forms described in subsection \ref{gluecohsub} does on
matching pairs of harmonic forms, i.e. on the space
\[ \harm^{m}_{y} = (\harm^{m}_{+,0} \times \harm^{m}_{-,0}) \cap\ycal{Z}^{m} .\]
We prove that any cohomology class on $M$ can be obtained
by gluing a matching pair of harmonic forms with a fixed gluing parameter $L$,
except when $L$ corresponds to an eigenvalue of a certain endomorphism that we
will define below.

\begin{thm}
\label{gluehodgethm}
Let $M_{+}$, $M_{-}$ have EAC metrics.
Considering $L$ as a parameter, the linear map
\begin{equation}
\label{lfixeq}
Y_{H} : \harm^{m}_{y} \to H^{m}(M), \;\:
(\psi_{+}, \psi_{-}) \mapsto Y_{H}(\psi_{+}, \psi_{-}, L)
\end{equation}
is an isomorphism except when $-2L$ is an eigenvalue of
\begin{equation}
\label{saprotoeq}
\pi_{E}(\partial_{+}^{-1}C_{+} + \partial_{-}^{-1}C_{-}) :
E^{m-1}_{\diag} \to E^{m-1}_{\diag} .
\end{equation}
\end{thm}

We can write $H^{m}(X)$ as an orthogonal direct sum
$A^{m}_{\gi} \oplus E^{m}_{\gi}$, where $A^{m}_{\gi}$ is the image of
$j_{\gi}^{*} : H^{m}(M_{\gi}) \to H^{m}(X)$.

Let $\diagcal{A}^{m} = \onecal{A}^{m} \cap \twocal{A}^{m}$. This is then
the space of harmonic representatives for~$A^{m}_{d}$. Similarly let
$\diagcal{E}^{m} = \onecal{E}^{m} \cap \twocal{E}^{m}$, and denote
the corresponding subspace of $H^{m}(X)$ by $E^{m}_{d}$.
Let $\pi_{E} : H^{m}(M) \to E^{m}_{\diag}$ denote the \ltwoorth{} projection.

Recall that $\partial_{\gi}$ denotes the boundary map in the long exact
sequence for relative cohomology~(\ref{relexacteq2}).
It is convenient to define an isomorphism
\[ C_{\gi} : E_{\gi}^{m-1} \to \im \partial_{\gi} \subseteq
H^{m}_{cpt}(M_{\gi})\]
as follows. For $\tau \in \gical{E}^{m-1}$ let $\psi$ be the unique element of
$\harm^{m}_{\gi,E}$ (the bounded exact harmonic forms on $M_{\gi}$)
such that $B_{\gi,e}(\psi) = \tau$.
If we take $\eta_{\pm}$ as defined in (\ref{etaeq}) and $\rho_{\gi}$ a cut-off
function for the cylinder of $M_{\gi}$ then
$\psi + d\eta_{\gi}(\psi) - d(\rho_{\gi}t_{\gi}\tau)$ has
compact support, so represents a class
$C_{\gi}([\tau]) \in H^{m}_{cpt}(M_{\gi})$.
This class is mapped to $0$ by $e_{\gi}$, so lies in the image of
$\partial_{\gi}$.
Composing $C_{\gi}$ with the inverse of
$\partial_{\gi} : E^{m-1}_{\gi} \to \im \partial_{\gi}$ gives an
endomorphism
$\partial_{\gi}^{-1} C_{\gi}$ of $E^{m-1}_{\gi}$.

\begin{rmk}
\label{pdlem}
$\partial_{\gi}^{-1}C_{\gi} : E^{m}_{\gi} \to E^{m}_{\gi}$ is self-adjoint,
and hence so is the endomorphism (\ref{saprotoeq}).
$C_{\gi}$ is independent of the choice of $\rho_{\gi}$, but depends on both
the metric and the cylindrical coordinate --
replacing $t_{\gi}$ by $t_{\gi} + \lambda$ adds $\lambda \partial_{\gi}$
to $C_{\gi}$.
\end{rmk}

\begin{proof}[Proof of theorem \ref{gluehodgethm}]
Consider the map
$(i^{*}_{+} \oplus i^{*}_{-}) : H^{m}(M) \to H^{m}(M_{+}) \oplus H^{m}(M_{-})$
in the Mayer-Vietoris sequence.
Recall that $L$ is fixed, so that $Y_{H}$ is regarded as a linear map
$\harm^{m}_{y} \to H^{m}(M)$. To show that it is an isomorphism
it suffices to show that $\im((i^{*}_{+} \oplus i^{*}_{-}) \circ Y_{H}) =
\im(i^{*}_{+} \oplus i^{*}_{-})$, and that
$Y_{H} : \ker((i^{*}_{+} \oplus i^{*}_{-}) \circ Y_{H}) \to
\ker (i^{*}_{+} \oplus i^{*}_{-})$ is an isomorphism.

$(i^{*}_{+} \oplus i^{*}_{-})Y_{H}(\psi_{+}, \psi_{-}) =
([\psi_{+}], [\psi_{-}])$ and it follows from the exactness of
the Mayer-Vietoris sequence that
$\im((i^{*}_{+} \oplus i^{*}_{-}) \circ Y_{H}) =
\im(i^{*}_{+} \oplus i^{*}_{-})$. It also follows that
$\ker ((i^{*}_{+} \oplus i^{*}_{-}) \circ Y_{H}) = \harm^{m}_{y,E}$, the
pairs of exact forms in $\harm^{m}_{y}$.

Thus the problem reduces to determining whether the restriction
\[ Y_{H} : \harm^{m}_{y,E} \to \ker(i^{*}_{+} \oplus i^{*}_{-}) \]
of (\ref{lfixeq}) is an isomorphism.
Given $\tau \in \diagcal{E}^{m-1}$ let $(\psi_{+},\psi_{-})$ be the unique
element of $\harm^{m}_{y,E}$ such that
$\tau = B_{+,e}(\psi_{+}) = -B_{-,e}(\psi_{-})$.
By the definition of $C_{\gi}$ and (\ref{yheq})
\begin{equation*}
Y_{H}(\psi_{+}, \psi_{-}) =
i_{+*}C_{+}([\tau]) + i_{-*}C_{-}([-\tau]) + 2L\delta([\tau]) .
\end{equation*}
Combining with (\ref{dieq})
\begin{equation}
\label{yharmeeq}
Y_{H}(\psi_{+}, \psi_{-}) =
\delta \left(\partial_{+}^{-1}C_{+}([\tau]) +
\partial_{-}^{-1}C_{-}([\tau]) + 2L[\tau] \right) .
\end{equation}
$\delta : H^{m-1}(X) \to H^{m}(M)$ is an isomorphism
$E^{m-1}_{\diag} \to \ker (i^{*}_{+} \oplus i^{*}_{-})$
and vanishes on the orthogonal complement of $E^{m-1}_{\diag}$. 
It follows that (\ref{yharmeeq}) gives an isomorphism
$\harm^{m}_{y,E} \to \ker (i^{*}_{+} \oplus i^{*}_{-})$ unless
$-2L$ is an eigenvalue of the endomorphism (\ref{saprotoeq}).
\end{proof}

\section{The gluing map}
\label{glueproofsec}

We will now make use of the topological results of the previous section to
study the gluing map for torsion-free \gtstr s. As in section
\ref{gluesetupsec} the setup is that $M_{+}$ and $M_{-}$ are EAC \gtmfd s
with a common cross-section $X$, and $M$ is their connected sum.
$\mathcal{M}$ denotes the moduli space of torsion-free \gtstr s on $M$, and
$G_{0}$ the space of gluing data.

In order to prove theorem \ref{maingluethm} we need to show that the gluing map
is invariant under $\ycal{D}$ (the identity component of the group of matching
pairs of EAC diffeomorphisms of $M_{+}$ and $M_{-}$) so that it is well-defined
on $\calg_{0} = G_{0}\ycal{D}/\ycal{D}$, show that $\calg_{0}$ is a smooth
manifold, and compute the derivative of the gluing map.

\subsection{Diffeomorphism invariance}

Note that the composition $\pi_{H} \circ Y : G_{0} \to H^{3}(M)$ of
the gluing map (\ref{upperyeq})
with the local diffeomorphism $\pi_{H} : \mathcal{M} \to H^{3}(M)$
is simply the restriction to $G_{0}$ of the map $Y_{H}$ given by
definition \ref{yhdefn}. We will use this first to show that $Y$
induces a well-defined map on the quotient $\mathcal{G}_{0}$.
Later we will determine the local
properties of $Y : \mathcal{G}_{0} \to \mathcal{M}$ from those of
$Y_{H} : (\ycal{X} \times \bbr)/\ycal{D} \to H^{3}(M)$.

\begin{prop}
\label{glueinvprop}
The map $Y : G_0 \to \mathcal{M}$ is $\ycal{D}$-invariant, so descends to a
well-defined continuous function
\begin{equation}
Y : \mathcal{G}_{0} \to \mathcal{M} .
\end{equation}
\end{prop}

\begin{proof}
We need to show that if $\phi \in \ycal{D}$ and
$(\varphi_{+}, \varphi_{-}, L) \in G_{0}$ such that
$\phi^{*}(\varphi_{+}, \varphi_{-}, L) \in G_{0}$ then
\[ Y(\varphi_{+}, \varphi_{-}, L) = Y(\phi^{*}(\varphi_{+}, \varphi_{-}, L)) .\]
The idea of the proof is to connect $(\varphi_{+}, \varphi_{-}, L)$ and
$\phi^{*}(\varphi_{+}, \varphi_{-}, L)$ by a path in~$G_{0}$. The image
under $Y$ of this path is the lift by the local diffeomorphism
$\pi_{H} : \mathcal{M} \to H^{3}(M)$ of a path in $H^{3}(M)$, which
is determined by propositions \ref{yhinvprop} and \ref{lderprop}.

Let $[0,1] \to \mathcal{M}$, $s \mapsto \phi_{s}$ be a path in $\ycal{D}$
connecting the identity to $\phi$, and take $k$ sufficiently large that
$\phi_{s}^{*}(\varphi_{+}, \varphi_{-}, L+k) \in G_{0}$ for all $s$.
By proposition \ref{yhinvprop} the path 
$s \mapsto Y(\phi_{s}^{*}(\varphi_{+}, \varphi_{-}, L+k)) \in \mathcal{M}$
is a lift of a constant path in $H^{3}(M)$, so
\[ Y(\phi^{*}(\varphi_{+}, \varphi_{-}, L+k)) =
Y((\varphi_{+}, \varphi_{-}, L+k)) . \]
By proposition \ref{lderprop} the paths
\begin{gather*}
s \mapsto Y(\phi^{*}(\varphi_{+}, \varphi_{-}, L+(1-s)k)) \in \mathcal{M} , \\
s \mapsto Y((\varphi_{+}, \varphi_{-}, L+(1-s)k)) \in \mathcal{M} .
\end{gather*}
are both lifts
of $s \mapsto Y_{H}(\varphi_{+}, \varphi_{-}, L+k) - 2ks\delta([\omega])$,
so in particular they have the same value at $s = 1$, which gives the result.
\end{proof}

\subsection{Deformations of EAC $G_{2}$-manifolds}
\label{eacdefsub}

In order to define coordinate charts for $\calg_{0}$ we first summarise
the deformation theory for EAC \gtmfd s developed in \cite[\S $6$]{jn1}.
Let $\calx_{\gi}$ be the space of torsion-free EAC \gtstr s on $M_{\gi}$
(with any exponential rate) and $\cald_{\gi}$ the group of EAC diffeomorphisms
isotopic to the identity.
Then the moduli space $\defstr_{\gi} = \calx_{\gi}/\cald_{\gi}$ is a smooth
manifold.

Any EAC torsion-free \gtstr{} $\varphi_{\gi}$ on $M_{\gi}$ is asymptotic to
$\Omega \pm dt_{\gi} \wedge \omega$, where $(\Omega, \omega)$ is a Calabi-Yau
structure on $X$.
This defines a natural boundary map
$B_{\gi} : \defstr_{\gi} \to \caln$, where $\caln$ is the moduli space of
Calabi-Yau structures on $X$.
Since $[\Omega] = j^{*}_{\gi}[\varphi_{\gi}]$ and
$\half[\omega]^{2} = j^{*}_{\gi}[*\varphi_{\gi}]$ it is clear that any
element in the image of the boundary map satisfies
\begin{equation}
\label{necceq}
[\Omega] \in A^{3}_\gi, \;\: [\omega]^{2} \in A^{4}_\gi ,
\end{equation}
where $A^m_\gi = \im(j^*_\gi : H^m(M_\gi) \to H^m(X))$ as before.
These conditions define a subset $\caln_{\gi,A} \subseteq \caln$.
The boundary map $B_{\gi}$ is a submersion onto its image, which is a
submanifold of $\caln$ and an open subset of $\caln_{\gi,A}$.

The proof of these results uses \emph{pre-moduli spaces} as coordinate charts.
There is a manifold $\gical{R}$ of torsion-free EAC \gtstr s near
$\varphi_{\gi}$, such that the natural map $\gical{R} \to \defstrgi$ is
a homeomorphism onto an open subset. The transition function between such maps
are smooth, so they can be used as coordinate charts.
Similarly there is a manifold $\calq$ of Calabi-Yau structures near
$(\Omega, \omega)$ such that $\calq \to \caln$ is a coordinate chart.

The subset $\calq_{\gi,A} \subseteq \calq$ defined by the equations
(\ref{necceq}) is a submanifold, and
the boundary map
\[ B_{\pm} : \gical{R} \to \calq \] 
is a submersion onto $\calq_{\gi,A}$.
Any tangent $(\sigma,\tau)$ at $(\Omega,\omega)$ to a curve of
$SU(3)$-structures satisfies the algebraic relations
\begin{gather}
\label{loneeq}
\sigma \wedge *\Omega - \tau \wedge \omega^{2} = 0 , \\
\label{ltwoeq}
\sigma \wedge \omega + \Omega \wedge \tau = 0 .
\end{gather}
The tangent space to $\calq$ consists of the harmonic tangents to the space
of $SU(3)$-structures, i.e.
\[ T_{(\Omega,\omega)}\calq = \harm_{SU} =
\{ (\sigma, \tau) \in \harm^{3}_{X} \times \harm^{2}_{X} :
\textrm{\eqref{loneeq} and \eqref{ltwoeq} hold} \} . \]
The tangent space to $\gical{R}$ consists of harmonic forms asymptotic
to elements of $\harm_{SU}$,
\[ T_{\varphi_{\gi}}\gical{R} = \harm^{3}_{\gi, cyl} =
\{ \psi \in \harm^{3}_{\gi,0} : B_{\gi}(\psi) \in \harm_{SU} \} . \]

\subsection{A coordinate chart}

Next we describe coordinate charts for $(\ycal{X} \times \bbr)/\ycal{D}$,
which contains $\mathcal{G}_{0}$ as an open subset.
$(\ycal{X} \times \bbr)/\ycal{D}$ is a principal $\bbr$-bundle over
$\calb = \ycal{X}/\ycal{D}$, so it suffices to show that $\calb$ is a manifold.

Let $\gical{M}$ be the moduli space of torsion-free EAC \gtstr s on $M_{\gi}$,
and $\mathcal{N}$ the moduli space of \cystr s on $X$. Let
$\ycal{M} \subseteq \onecal{M} \times \twocal{M}$ be the pairs of
diffeomorphism classes of EAC torsion-free \gtstr s whose boundary images in
$\mathcal{N}$ match.

\begin{prop}
\label{coverprop}
Let $M$ be the gluing of a pair of EAC \gtmfd s $M_{\gi}$. If $b^{1}(M) = 0$
then the natural projection
\begin{equation}
\label{glueprojeq}
\calb \to \ycal{M},
\;\: (\varphi_{+}, \varphi_{-})\ycal{D} \mapsto
(\varphi_{+}\onecal{D}, \varphi_{-}\twocal{D})
\end{equation}
is a local homeomorphism.
\end{prop}

\noindent
First, we find charts for $\ycal{M}$.

\begin{prop}
\label{mymfdprop}
$\ycal{M}$ is a submanifold of $\onecal{M} \times \twocal{M}$.
\end{prop}

Each point in $\ycal{M}$ can be represented by
a matching pair of torsion-free \gtstr s $(\varphi_{+}, \varphi_{-})$,
asymptotic to a Calabi-Yau structure $(\Omega,\omega)$ on $X$.
Let $\gical{R}$ be the pre-moduli space of torsion-free EAC \gtstr s near
$\varphi_{\gi}$. 

\begin{defn}
The pre-moduli space of matching pairs of torsion-free EAC \gtstr s near
$(\varphi_{+},\varphi_{-})$ is a neighbourhood $\ycal{R}$
of $(\varphi_{+},\varphi_{-})$ in
$\ycal{X} \cap (\onecal{R} \times \twocal{R})$.
\end{defn}

To use $\ycal{R}$ as a coordinate chart we first need to show that its
image under the boundary map is a manifold.
The intersection
\[ \mathcal{Q}_{\diag,A} = \mathcal{Q}_{+,A} \cap \mathcal{Q}_{-,A} \]
consists of $(\Omega', \omega') \in \mathcal{Q}$ such that
$[\Omega'] \in A^{3}_{\diag}, [\omega'^{2}] \in A^{4}_{\diag}$.

\begin{lem}
\label{qdamfdlem}
$\mathcal{Q}_{\diag,A} \subseteq \mathcal{Q}$ is a submanifold.
\end{lem}

\begin{proof}
The proof of proposition \cite[Proposition $6.2$]{jn1}, which states
that each of $\calq_{\gi,A}$ is a manifold, can be recycled.
\end{proof}

\begin{proof}[Proof of proposition \ref{mymfdprop}.]
The group $\xcal{D}$ of diffeomorphisms of $X$ isotopic to the identity
acts trivially on $\mathcal{Q}$,
so for $(\psi_{+}, \psi_{-}) \in \onecal{R} \times \twocal{R}$
\[ B_{+}(\psi_{+}) \; \xcal{D}\textrm{-equivalent to } B_{-}(\psi_{-})
\iff B_{+}(\psi_{+}) = B_{-}(\psi_{-}) . \]
Hence  $\ycal{R}$ is homeomorphic to a neighbourhood of
$(\varphi_{+}\onecal{D}, \varphi_{-}\twocal{D})$ in $\ycal{M}$,
and it suffices to prove that
$\ycal{R}$ is a submanifold of $\onecal{R} \times \twocal{R}$.

By lemma \ref{qdamfdlem} the image of $\mathcal{Q}_{\diag,A}$ in
$\mathcal{Q}_{+,A} \times \mathcal{Q}_{-,A}$ under the diagonal map
is a submanifold.
$\ycal{R} \subseteq \onecal{R} \times \twocal{R}$ is the inverse image
of $\mathcal{Q}_{\diag,A} \subseteq \mathcal{Q}_{+,A} \times \mathcal{Q}_{-,A}$
under the submersion $B_{+} \times B_{-} : \onecal{R} \times \twocal{R} \to
\mathcal{Q}_{+,A} \times \mathcal{Q}_{-,A}$,
so it is a submanifold.
\end{proof}

\begin{proof}[Proof of proposition \ref{coverprop}.]
Let $(\varphi_{+},\varphi_{-}) \in \ycal{X}$, and
$\ycal{R}$ the pre-moduli space of nearby matching pairs.
Because $\ycal{R}$ is a coordinate chart for $\ycal{M}$, any element 
of $\ycal{X}$ near $(\varphi_{+},\varphi_{-})$ can be written as 
$(\phi_{+}^{*}\psi_{+}, \phi_{-}^{*}\psi_{-})$, with
$(\psi_{+}, \psi_{-}) \in \ycal{R}$ and $\phi_{\gi} \in \gical{D}$ close
to $id$.

Let $Aut_{0}(X) \subset \xcal{D}$ be the identity component of the subgroup of
automorphisms of the Calabi-Yau manifold $X$
(this is actually independent of $(\psi_{+}, \psi_{-}) \in \ycal{R}$, 
cf. \cite[Proposition $4.5$]{jn1}).
The matching condition for $(\phi_{+}^{*}\psi_{+}, \phi_{-}^{*}\psi_{-})$
implies that $B(\phi_{-})^{-1}B(\phi_{+}) \in Aut_{0}(X)$, where
$B(\phi_{\gi})$ denotes the asymptotic limit of $\phi_{\gi}$.

$Aut_{0}(X)$ is a closed subgroup of the isometry group
of $X$, so it is compact (see Myers and Steenrod \cite{myers39}).
Because $X$ is Ricci-flat the Lie algebra
of $Aut_{0}(X)$ corresponds to the space $\harm^{1}_{X}$ of harmonic $1$-forms
on $X$. Because $X$ is Ricci-flat these are parallel, so the group is abelian.
Similarly, the Lie algebras of the automorphism groups $Aut_{0}(M_{\gi})$ of
the EAC \gtmfd s $M_{\gi}$ (which are independent of
$\varphi_\gi \in \calr_\gi$) correspond to the bounded harmonic
\mbox{$1$-forms} $\harm^{1}_{\gi,0}$.
The image of $\harm^{1}_{\gi,0}$ under the boundary map $B_{\gi}$ is the
space of harmonic representatives of
$A^{1}_{\gi} = \im (j_{\gi}^{*} : H^{1}(M_{\gi}) \to H^{1}(X))$.
Each of these is a half-dimensional subspace of $H^{1}(X)$ according to
\cite[Proposition $5.15$]{jn1}, while their intersection is
$A^{1}_{d} = \im (j^{*} : H^{1}(M) \to H^{1}(X))$. As we assume $b^{1}(M) = 0$
it follows that $H^{1}(X) = A^{1}_{+} \oplus A^{1}_{-}$.
Hence $Aut_{0}(X)$ is generated by the images $B(Aut_{0}(M_{\gi}))$.

\hfuzz=2pt
It follows that $B(\phi_{-})^{-1}B(\phi_{+}) = B(\phi'_{-})^{-1}B(\phi'_{+})$
for some $\phi'_{\gi} \in Aut_{0}(M_{\gi})$. Then
$(\phi'_{+}\phi_{+}, \phi'_{-}\phi_{-})$ is a matching pair of diffeomorphisms,
\hfuzz=0.1pt
so $(\phi_{+}^{*}\psi_{+}, \phi_{-}^{*}\psi_{-})$ is $\ycal{D}$-equivalent
to $(\psi_{+}, \psi_{-})$. Hence the image of $\ycal{R}$ is open in $\calb$,
so $\ycal{R}$ can be used as a coordinate chart for $\calb$ too.
\end{proof}

\noindent
Because $\ycal{R}$ is a coordinate chart for $\calb$ when $b^{1}(M)= 0$,
the natural maps
\begin{equation}
\label{gluecharteq}
\ycal{R} \times \bbr \to (\ycal{X} \times \bbr)/\ycal{D}
\end{equation}
can then be used as local trivialisations for $(\ycal{X} \times \bbr)/\ycal{D}$
as a principal \mbox{$\bbr$-bundle}.

\begin{rmk}
It is possible to show that $\calb \to \ycal{M}$ is a covering map
when $b^{1}(M) = 0$.
In general the connected components of the fibres are isomorphic to $H^{1}(M)$.
\end{rmk}

\subsection{The derivative of the gluing map}

Since $\pi_{H} : \mathcal{M} \to H^{3}(M)$ is a local diffeomorphism
the local behaviour of the gluing map $Y : \mathcal{G} \to \mathcal{M}$
is determined by that of $Y_{H} = \pi_{H} \circ Y$. $Y_{H}$ is just
the gluing map for cohomology from definition \ref{yhdefn}, so can be defined
on all of $(\ycal{X} \times \bbr)/\ycal{D}$. We compute the
derivative.

\begin{prop}
\label{yhderprop}
Given $(\varphi_{+}, \varphi_{-}) \in \ycal{X}$ the derivative of
\[ Y_{H} : (\ycal{X} \times \bbr)/\ycal{D} \to H^{3}(M) \]
at $(\varphi_{+}, \varphi_{-}, L)$ is bijective
for all sufficiently large values of $L$.
\end{prop}

\begin{proof}
We make the simplifying assumption that $b^{1}(M) = 0$. The tangent space
to the pre-moduli space of matching torsion-free \gtstr s is
\[ T_{(\varphi_{+}, \varphi_{-})}\ycal{R} = \harm^{m}_{y,cyl} , \]
the space of matching harmonic forms $(\psi_{+}, \psi_{-})$ whose common
boundary value $B(\psi)$ lies in $\harm_{SU}$.
The condition that $b^{1}(M) = 0$ implies that the common boundary value
of any $(\psi_{+}, \psi_{-}) \in \harm^{m}_{y}$ automatically satisfies
(\ref{loneeq}).  $\harm^{m}_{y,cyl} \subset \harm^{m}_{y}$ therefore has
codimension $1$, and we can take $\{ (\psi_{+}, \psi_{-}) \in \harm^{m}_{y,E}:$
$B(\psi) \in \bbr[\omega] \}$ as a direct complement.
Let
\begin{equation*}
Y'_{H} : \ycal{R} \times \bbr \to H^{3}(M)
\end{equation*}
be the representation of $Y_{H}$ in the coordinate chart (\ref{gluecharteq}),
and consider the derivative
\[ (DY'_{H})_{(\varphi_{+},\varphi_{-},L)} :
\harm^{3}_{y,cyl} \times \bbr \to H^{3}(M) . \]
The restriction of $DY'_{H}$ to $\harm^{3}_{y,cyl} \times 0$ is just
$Y_{H}$ (\ref{lfixeq}),
while on $0 \times \bbr$ it is $h \mapsto 2h\delta([\omega])$.
By a slight modification of the proof of theorem \ref{gluehodgethm} we find
that 
\[ (i^{*}_{+} \oplus i^{*}_{-}) \circ Y_{H} : \harm^{3}_{y,cyl}
\to \im(i^{*}_{+} \oplus i^{*}_{-}) \]
is surjective with kernel
$\harm^{m}_{y,cyl,E} = \harm^{m}_{y,cyl} \cap \harm^{m}_{y,E}$,
and that if we identify
$\harm^{m}_{y,cyl,E} \times \bbr \leftrightarrow E^{2}_{d}$ by
$(\psi_{+}, \psi_{-},h) \mapsto [B_{e}(\psi) + \frac{h}{L}\omega]$ then
$DY'_{H} : \harm^{m}_{y,cyl,E} \times \bbr \to
\ker (i^{*}_{+} \oplus i^{*}_{-})$ is identified with
\[ E^{2}_{d} \to \ker (i^{*}_{+} \oplus i^{*}_{-}) ,
\tau \mapsto \delta(2L\tau + F(\tau)) \]
for some endomorphism $F$ of $E^{2}_{d}$
($F$ is the composition of (\ref{saprotoeq}) with the projection to the
orthogonal complement of $\bbr[\omega]$ in $E^{2}_{d}$).
Hence $DY'_{H}$ is an isomorphism except when $-2L$ is an eigenvalue of~$F$.
\end{proof}

We can define $G \subseteq G_{0}$ to be the subset of gluing data
$(\varphi_{+}, \varphi_{-}, L)$ for which the gluing parameter $L$ is
sufficiently large to ensure invertibility of the derivative of the gluing map.
The quotient $\mathcal{G} = G\ycal{D}/\ycal{D}$ is an open subset of
$\mathcal{G}_{0}$, and $Y : \mathcal{G} \to \mathcal{M}$ is a local
diffeomorphism.
This completes the proof of theorem \ref{maingluethm}.

\section{Boundary points of the moduli space}
\label{glueboundarysec}

To conclude we describe how to attach boundary points to the
moduli space of torsion-free \gtstr s of a compact \gtmfd{} $M$ obtained by
gluing, outlining a proof of theorem \ref{gluebdrythm}.

Let the compact manifold $M$ be the gluing of two EAC \gtmfd s $M_{\gi}$
as before, and
\[ Y : \calg \to \defstr \]
the gluing map for torsion-free \gtstr s.
The gluing space $\calg$ is a fibre bundle over $\calb$ with typical
fibre~$\bbrp$. It can be considered as the interior of a
topological manifold $\calgbar$ with boundary
$\calb$ `at infinity' by adding a limit point to each of the fibres.
We aim to add a boundary to $\defstr$ so that $Y$ extends to a local
homeomorphism $Y : \calgbar \to \defstrbar$.

Assume that $b^{1}(M) = 0$, and let $\ycal{R}$ be the pre-moduli space of
matching pairs of torsion-free EAC
\gtstr s near $(\varphi_{+}, \varphi_{-}) \in \ycal{X}$.
We can interpret the proof of proposition \ref{yhderprop} as stating that
$(i^{*}_{+} \oplus i^{*}_{-}) \circ Y_{H} : \ycal{R} \times \bbr
\to \im(i^{*}_{+} \oplus i^{*}_{-})$ is a submersion, and that
\begin{equation}
\label{coneeq}
Y_{H} : \ycal{R}' \times (L_{1}, \infty) \to K
\end{equation}
is a local diffeomorphism for $L$ sufficiently large, where
\[ \ycal{R}' = \{ (\varphi'_{+}, \varphi'_{-}) \in \ycal{R} :
i^{*}_{\gi}[\varphi'_{\gi}] = i^{*}_{\gi}[\varphi_{\gi}] \} \] 
and $K = \{ [\alpha] \in H^{3}(M) :
i^{*}_{\gi}[\alpha]= i^{*}_{\gi}[\varphi_{\gi}] \}$
(an affine translate of $\delta(H^{2}(X)) \subseteq H^{3}(M)$).
But we can make a stronger statement. The map (\ref{coneeq}) has the form
\begin{equation*}
Y_H(\varphi'_{+}, \varphi'_{-}, L) =
Y_H(\varphi'_{+}, \varphi'_{-}, 0) + 2L\delta([\omega']) ,
\end{equation*}
where $\omega'$ is the Kähler form of the common boundary value of
$(\varphi'_{+},\varphi'_{-}) \in \ycal{R}'$.
The second term maps out an open cone in $\delta(H^{2}(X))$.
For large enough $L_{1}$ the second term dominates, and (\ref{coneeq}) is
a diffeomorphism onto approximately an open affine cone in $K$. Hence
\begin{equation}
\label{bdrycharteq}
Y : \ycal{R} \times (L_{1}, \infty) \to \defstr
\end{equation}
is not just a local diffeomorphism, but a diffeomorphism onto its image
for large $L_{1}$.
Since $\ycal{R}$ are coordinate charts for $\calb$,
one could try to use (\ref{bdrycharteq}) as coordinate charts to make
$\defstr \cup \calb$ a manifold with boundary. The problem is that the
resulting topology need not be Hausdorff; different points of $\calb$ could
\emph{a priori} arise as the limit of the same path in $\defstr$.
This difficulty can be resolved by proving that the property of `defining
the same boundary point' is an equivalence relation on $\calb$ and that the
quotient $\hat\calb$ is covered by~$\calb$. Then $\hat\calb$ is a manifold,
and one can use $\defstrbar = \defstr \cup \hat\calb$ in the statement
of theorem \ref{gluebdrythm}.

This outline can be expanded to a full proof of theorem \ref{gluebdrythm}.
The details can be found in \cite[\S 6.4]{jnthesis}, 
but are not included here as they amount to a rather tedious inspection of
the charts (\ref{bdrycharteq}).

\begin{rmk}
We could give the topological manifolds $\calgbar$ and $\defstrbar$ smooth
structures by choosing an identification of $(0, \infty]$ with
a half-open interval $[0, 1)$, but it is not clear if there is a natural
choice of smooth structure.
\end{rmk}

\emph{Acknowledgements.}
I thank Alexei Kovalev for helpful discussions, Jonny Evans for constructive
comments, and the Swedish Royal Academy of Sciences funds for financial
support.

\bibliographystyle{plain}
\bibliography{g2geom}

\end{document}